\documentclass[onefignum,onetabnum]{siamart190516}
\usepackage[utf8]{inputenc}
\usepackage{subfigure,mathtools,mathrsfs} 
\usepackage[english]{babel}
\usepackage[font=small]{caption}
\newcommand{\rone}[1]{#1}
\newcommand{\rtwo}[1]{#1}
\newcommand{\rthr}[1]{#1}
\newcommand{\norm}[1]{\left\lVert #1 \right\rVert} 

\newcommand{\bm}[1]{\boldsymbol{#1}}    

\usepackage{lipsum}
\usepackage{amsfonts}
\usepackage{graphicx}
\usepackage{epstopdf}

\headers{Hybrid Stokes Drag-IB Method for Slender Bodies}{Ondrej Maxian and Charles S. Peskin}

\title{An immersed boundary method with subgrid resolution and improved numerical stability applied to slender bodies in Stokes flow}

\author{Ondrej Maxian\thanks{Courant Institute of Mathematical Sciences, New York University, 251 Mercer St., New York, NY 10012
  (\email{om759@nyu.edu}).}
\and Charles S. Peskin\thanks{Courant Institute of Mathematical Sciences,  New York University, 251 Mercer St., New York, NY 10012
  (\email{peskin@cims.nyu.edu}).}}

\usepackage{amsopn}

\ifpdf
\hypersetup{
  pdftitle={An immersed boundary method with subgrid resolution and improved numerical stability applied to slender bodies in Stokes flow},
  pdfauthor={Ondrej Maxian and Charles S. Peskin}
}
\fi

\hypersetup{
colorlinks = true, linkcolor = blue, citecolor = black, urlcolor=black}

\begin{document}

\maketitle

\begin{abstract}
  The immersed boundary method is a numerical and mathematical formulation for solving fluid-structure interaction problems. It relies on solving fluid equations on an Eulerian fluid grid and interpolating the resulting velocity back onto immersed structures. To resolve slender fibers, the grid spacing must be on the order of the fiber radius, and thus the number of required grid points along the filament must be of the same order as the aspect ratio. Simulations of slender bodies using the IB method can therefore be intractable. A technique is presented to address this problem in the context of Stokes flow. The velocity of the structure is split into a component coming from the underlying fluid grid, which is coarser than normally required, and a component proportional to the force (a drag term). The drag coefficient is set so that a single sphere is represented exactly on a grid of arbitrary meshwidth. Implicit treatment of the drag term removes some of the stability restrictions normally associated with the IB method. This comes at a loss of accuracy, although tests are conducted that show 1-2 digits of relative accuracy can be obtained on coarser grids. After its accuracy and stability are tested, the method is applied to two real world examples: fibers in shear flow and a suspension of fibers. These examples show that the method can reproduce existing results and make reasonable predictions about the viscosity of an aligned fiber suspension.
\end{abstract}



\section{Introduction}
\label{sec:intro}
Many problems from biology are fundamentally interaction problems between bodies and the fluid in which they are immersed. The immersed boundary method \cite{peskin2002acta} offers a mathematical and numerical approach to treat such problems. Although first introduced to study flow patterns around heart valves \cite{peskin1972flow}, the IB method has also been applied to systems at larger and smaller length scales, such as jellyfish swimming \cite{jellyfish} and cellular motility \cite{wandableb13}, respectively. The IB method relies on a set of discrete interaction equations that communicate the force from the immersed body to a volumetric fluid grid. Fluid equations are solved on the grid, and the resulting velocity is interpolated back to the immersed structure.

Applying the IB method to multiscale problems can be prohibitively expensive because of the need to solve fluid equations in the entire volume. One example of this is the simulation of slender actin and microtubular filaments that comprise the cellular cytoskeleton \cite{alberts}. In previous work on slender body IB simulations, Bringley and Peskin \cite{tbring08}, found that every immersed point (``marker'') in the IB method has an associated radius, now generally referred to as the ``hydrodynamic radius'' \cite{rigidmblob}. This radius is on the order of the grid size. Thus to resolve the radius of a slender filament, an extremely fine (relative to the length of the filament) grid must be used for the fluid domain. This in turn leads to an increased number of markers for the filament (in comparison to the number that would be used on a coarser fluid grid), since the two grid spacings need to be comparable \cite{tbring08}. Making the resolution finer on an immersed elastic structure increases the numerical stiffness, especially if the structure resists bending so that fourth derivatives are involved \cite{leveque2007finite}. Thus in order to do large-scale simulations on slender bodies with an IB method, new numerical methods are needed that can give better accuracy on coarser grids and/or overcome the numerical stiffness problem.

One obvious way to address the spatial resolution issue is to employ adaptive mesh refinement. Griffith et al. \cite{griffith2007adaptive} developed an adaptive, distributed memory parallel implementation of the IB method (IBAMR), which has been used to model slender semiflexible diatom chains \cite{ngfauci} and elastic rings \cite{griffith2012ring}. In the limit of many fibers, however, adaptive mesh refinement breaks down as the fine mesh occupies most of space. Recently, Wiens and Stockie developed another distributed memory implementation of the IB method that employs a uniform mesh. By relaxing the incompressibility constraint and using a pseudo-compressible fluid solver for the Navier-Stokes equations \cite{wiens2015efficient}, they obtain a series of smaller, tridiagonal, linear systems to be solved at every timestep and therefore weak scalability. Because of this, they were able to simulate a suspension of 256 intrinsically curved fibers \rthr{with near optimal parallel scaling} \cite{wiens2015simulating}.

Implicit IB methods offer another way to reduce the cost of the IB method for slender bodies. Typically, implicit schemes are application-dependent and are based on solving large linear systems for the configuration of the immersed structure. Some existing schemes are based primarily on using semi-implicit solvers for the Navier-Stokes equations \cite{lee2010immersed, hou2008efficient} and time-lagging of the spreading and interpolation operators \cite{mori2008implicit, ceniceros2009efficient}. The latter approach greatly simplifies the design of implicit schemes without changing stability considerations \cite{fogguy07}. Approaches to accelerate the solution of large mixed Lagrangian-Eulerian linear systems that arise in these schemes include improved Jacobian-free Krylov solvers \cite{le2009implicit} and geometric multigrid as a preconditioner for a Krylov solver \cite{guy2015geometric}. 

In the present paper, we consider the case in which small length scales imply that the Stokes equations govern the fluid flow. A completely different approach available in this context is to eliminate the fluid grid entirely and use a Green's function for the steady state Stokes equations. For slender filaments, slender body theory \rone{(SBT)} \cite{keller1976slender, johnson80} gives an equation for the evolution of the fiber center line that is asymptotic in the aspect ratio $\epsilon$. This reduces the complexity of the solver from 3D to 1D. However, \rthr{the asymptotics of SBT break down} close to the filament and near the filament ends \cite{mori2018theoretical}. The non-local aspects of slender body theory are challenging in that they involve $\mathcal{O}(N^2)$ interaction terms (which must be treated with fast multipole methods or Ewald splitting for linear algorithm time) and singular integrals along fiber centerlines (which must be treated with appropriate quadrature schemes or regularization) \cite{ts04, ehssan17}. 

\rtwo{Another Green's function based approach to the problem is to approximate a fiber by a distribution of spheres and use analytical formulas for the mobility of a sphere in Stokes flow. For example, the force coupling method (FCM) represents each sphere by a Gaussian blob of force, where the width of the Gaussian $\sigma=r/\sqrt{\pi}$ is chosen to give the correct flow for a sphere of radius $r$ \cite{fcm03, fcm10}, and the Rotne-Prager-Yamakawa mobility (RPY) tensor treats the sphere as a surface delta source with radius $r$ \cite{rpyOG, donRPY}. For periodic systems, both of these methods require some grid-based solver which must be refined as the particle size decreases. For example, in periodic simulations that use FCM, the grid size is again recommended to be on the order of the fiber radius, and this can be prohibitive for systems of many slender fibers \cite{fcm10}.

The most rigorous way to coarsen the grid while controlling accuracy in this context is via Ewald splitting and the spectral Ewald method \cite{sham17Ew, donRPY}. Here the action of the Green's function is split into a smooth ``far field'' which rapidly decays in Fourier space (handled on a coarse grid via performing a non-uniform discrete Fourier transform), and a non-smooth ``near field'' which rapidly decays in real space \cite{lind11Ew, ewaldOG}. The method can be made log-linear in the number of particles by choosing a near field that is nonzero for $\mathcal{O}(1)$ neighbors. This makes the far field less smooth as the number of particles increases, which means the grid needed to compute it must be refined. }

Our goal here is to develop an immersed boundary method for slender bodies in Stokes flow that allows for a coarser spatial grid, decreases the overall temporal stiffness of the method, and achieves subgrid resolution. \rtwo{At the same time, we would like the method to be linear in the number of particles regardless of the chosen grid size}. In Section 2, we show how this is possible by combining the immersed boundary velocity from the grid with a Stokes drag term so that the total drag on a marker is exact for any given physical radius (not determined by the grid). In Section 3, we verify the accuracy of the method by computing the drag on an ellipsoid of large aspect ratio. We represent the ellipsoid by a 1D array of markers with non-constant radius (in accordance with its true geometry) and compare the result we obtain to slender body theory. In Section 4, we show that introducing the drag term also decreases the overall numerical stiffness of the IB method, and that the linear system to be solved depends only on the \textit{Lagrangian} quantities. Finally, we verify our method in Sections 5 and 6 by simulating a single fiber in a shear flow, and then measuring the viscosity of a suspension of up to 640 aligned fibers. 

\section{Method formulation}
\rone{In this section, we lay out our new IB method for slender bodies. This necessarily begins with a description of the classical IB method, followed by a discussion of how we can leverage the Stokes drag formula for spheres to obtain the velocity of a slender body from a coarse grid. We conclude this section by discussing the different possible temporal integrators for our IB method. }

\subsection{IB Method}
We consider the immersed boundary method in the special case in which the Stokes equations govern the fluid mechanics. That is, 
\begin{gather}
\label{eq:Stokeseqnsone}
\mu \Delta \bm{u} - \nabla p + \bm{f} = \bm{0}, \\
\label{eq:Stokeseqnstwo}
\nabla \cdot \bm{u} = 0,
\end{gather}
where $\mu$ is the fluid viscosity, $\bm{u}$ is the fluid velocity, $p$ is the pressure, and $\bm{f}$ is the external forcing. In an IB-type formulation, the immersed structures (in our case slender fibers) are represented using a Lagrangian description, while the fluid quantities $\bm{u}$ and $p$ are Eulerian variables (i.e. they are functions of fixed Cartesian coordinates $\bm{x}$ and time $t$). Throughout this manuscript, we follow the typical IB convention in referring to Lagrangian quantities using capital letters and Eulerian quantities using lowercase letters. 

Let $\Gamma$ denote the Lagrangian domain, parameterized by material parameter $q$, and let $\Omega$ denote the fluid domain (a uniform grid with spacing $h$). Suppose that the Lagrangian configuration has force density $\bm{F}(q)$. Then the first interaction equation for the force density is 
\begin{equation}
\label{eq:spread}
\bm{f}(\bm{x}) = \int_{\Gamma}  \bm{F}(q,t) \delta_h (\bm{x} - \bm{X}(q,t))\, dq  \approx \sum_{q} \bm{F}(q,t) \delta_h (\bm{x} - \bm{X}(q,t)) \Delta q = \bm{\mathcal{S}}\left(\bm{X}\right) \bm{F} . 
\end{equation}
Eq.\ \eqref{eq:spread} defines the \textit{spreading} operator $\bm{\mathcal{S}}$ and is used to calculate the force density $\bm{f}$ on the grid from the Lagrangian force density on the structure. $\delta_h$ is the classical 4 point regularized delta function of \cite{peskin2002acta}. 
After spreading the force to the grid via $\bm{f}=\bm{\mathcal{S}}\bm{F}$, the Stokes equations are solved on a periodic domain using a Fourier-spectral fluid solver (not a finite difference method) \cite{specflow}. The structure is updated with the local fluid velocity, satisfying the no-slip boundary condition,
\begin{align}
\label{eq:interp}
\frac{d \bm{X}}{dt} &= \bm{U}(\bm{X}(q,t)) =\bm{\mathcal{S}}^* \left(\bm{X}\right)\bm{u}\\[2 pt] \nonumber & =  \sum_{\bm{x}} \bm{u}(\bm{x},t) \delta_h (\bm{x} - \bm{X}(q,t)) h^3 \approx \int_{\Omega} \bm{u} (\bm{x},t) \delta_h(\bm{x}-\bm{X}(q,t))\, d\bm{x}. 
\end{align}
In Eq.\ \eqref{eq:interp}, we have defined the \textit{interpolation} operator $\bm{\mathcal{S}}^*$ that acts on the fluid grid velocity $\bm{u}$ to obtain the velocity at the structure points $\bm{U}$. It is important for energy conservation that the interpolation operator be the adjoint of the spreading operator \cite{peskin2002acta}. 

\subsection{Modifications for slender bodies}

As discussed in Section \ref{sec:intro}, simulations of slender bodies with the IB method are often intractable because of the multiple length scales in the problem. To resolve the diameter of the body, a fine grid must be used relative to the length of the body. Once the grid size becomes smaller, the spacing between markers also decreases, and so does the timestep necessary for stability. 

Our goal here is to leverage analytical results to obtain an implicit immersed boundary method with subgrid resolution. Since the fluid is governed by the Stokes equations, we know analytically the Stokes drag law for a spherical body
\begin{equation}
\label{eq:sdrag}
\bm{U} = \frac{\bm{F}}{6\pi \mu a}, 
\end{equation}
where $\bm{U}$ is the velocity of the body induced by a force $\bm{F}$, $\mu$ is the fluid viscosity, and $a$ is the radius of the body. While this paper is restricted to Stokes flow, drag coefficients for nonzero Reynolds number have been estimated empirically and could be used in principle in place of Eq.\ \eqref{eq:sdrag} (see \cite{clift2005bubbles} for empirical expressions). 

Previous studies of the IB method \cite{tbring08, rigidmblob} have shown that every immersed boundary marker has an associated numerical radius, $R_h$. Bringley and Peskin \cite{tbring08} systematically measured this quantity on an unbounded grid using several different spreading delta functions and fluid solvers (spectral vs. finite difference). They found that, for the standard 4 point delta function the mean $R_h \approx 1.2-1.3h$, with larger values for a spectral fluid solver than a finite difference one. In addition, they determined that the optimal spacing for an array of markers that represent a cylinder is $R_h$. We will also measure this quantity in Appendix A, with the result that $R_h=1.33h$ (this is the value we use throughout this paper). 

In a traditional IB formulation of the immersed fiber problem, the grid size is set so that the \textit{physical radius matches the hydrodynamic radius}. That is, on the particular grid chosen, Eq.\ \eqref{eq:sdrag} holds for a single point with $a=R_h=R$, where $R$ is the physical radius of the point. For slender bodies, using this technique can often render computations intractable. For example, actin filaments are typically about 1 $\mu$m long and have a radius $ \leq 5$ nm. Using the traditional IB formulation \rone{(which assumes no-slip conditions hold along the surface of the actin filament, see \cite{brem91} for contrary views),} this translates to 260 grid points and $\mathcal{O}(100)$ markers, just along the length of one filament. Thus simply simulating a 3D box that can hold 2 filaments aligned in each coordinate direction requires a $512^3$ grid, making detailed parameter studies impossible. 

We take a different approach. Suppose that we want to represent exactly a sphere with radius $R$ on a grid of spacing $h$, where $R$ is not necessarily equal to $R_h$. In fact, we want the physical radius, $R$, to be less than the grid hydrodynamic radius, $R_h$, so that we can solve the problem on a coarser grid. We seek to use the Stokes drag law, Eq.\ \eqref{eq:sdrag}, to correct the velocity from the immersed boundary method. More formally, if we apply a force of strength $\bm{F}$, Eq.\ \eqref{eq:sdrag} should hold with $a=R$. We can accomplish this by splitting the velocity into an IB part and a Stokes drag part, 
\begin{equation}
\label{eq:mobeq}
\bm{U} = \frac{\bm{F}}{6\pi \mu R} = \frac{\bm{F}}{6\pi \mu R_h} + \frac{\bm{F}}{6\pi \mu R_c},
\end{equation}
where $R_c$ is a radius that we use to correct the velocity of the sphere that comes from the IB method. We can solve Eq.\ \eqref{eq:mobeq} for the correction radius $R_c$, 
\begin{gather}
\frac{1}{R}=\frac{1}{R_h}+\frac{1}{R_c},\\[4 pt]
\label{eq:Rs}
R_c = \frac{R_hR}{R_h-R}. 
\end{gather}
The idea is therefore as follows. Begin with a grid spacing $h$ and therefore a hydrodynamic radius $R_h$. Then use Eq.\ \eqref{eq:Rs} to compute the required correction radius so that the sphere moves as if it has radius $R$. This involves replacing the IB velocity formula, which is traditionally $\bm{U} = \bm{\mathcal{S}}^* \bm{u}$, Eq.\ \eqref{eq:interp}, with
\begin{gather}
\label{eq:hybridupdate}
\bm{U} = \bm{\mathcal{S}}^* \bm{u} + \xi \bm{F},\\[4 pt]
\label{eq:dragc}
\xi = \frac{1}{6\pi \mu R_c}. 
\end{gather}
We make the following observations about Eq.\ \eqref{eq:Rs}. First, if $R_h=R$, then $R_c=\infty$ and it naturally follows from Eq.\ \eqref{eq:dragc} that $\xi=0$. This means there is no correction (and we are using the traditional IB method) when the fiber is fully resolved. On the other hand, when $R_h \rightarrow \infty$ (the fluid grid is entirely removed), we get $R_c=R$, and the Stokes drag model is recovered.

\subsection{Temporal discretization}
As we show in Section \ref{sec:stability}, the introduction of the drag term in Eq.\ \eqref{eq:hybridupdate} on a fixed grid \textit{increases} the overall numerical stiffness of the problem. For flexible fibers, this is because the force generally depends on fourth derivatives of $\bm{X}$. Eq.\ \eqref{eq:hybridupdate} lends itself to implicit treatment, however, most especially for forcing that is a linear function of position (e.g. bending forces). We can temporally discretize Eq.\ \eqref{eq:hybridupdate} as
\begin{equation}
\label{eq:gendisc}
\frac{\bm{X}^{n+1}-\bm{X}^n}{\Delta t} = \bm{\mathcal{S}}^*(\bm{X}^n) \bm{u}^n+ \xi \bm{F}(\bm{X}^{n+1}),
\end{equation}
which, depending on the form of $\bm{F}$, might be a linear system that is easily solvable for $\bm{X}^{n+1}$. In the subsequent examples, we will consider forcing that is both a linear function (bending resistance) and nonlinear function (fiber tension) of the position to illustrate when this method is most applicable. In the bending resistance case, a sparse banded linear system emerges, which can be solved in $\mathcal{O}(N)$ operations, where $N$ is the number of immersed points. For fiber tension, force is a nonlinear function of $\bm{X}$, and Newton's method must be used. Each Newton step, however, requires only the solution of a  banded linear system. 

In either case, we note that the non-local hydrodynamic interactions in Eq.\ \eqref{eq:gendisc} are treated \textit{explicitly}. Thus, in the case of multiple fibers, Eq.\ \eqref{eq:gendisc} can be solved fiber by fiber. The complexity of solving Eq.\ \eqref{eq:gendisc} therefore scales linearly with the number of points (assuming each solve is banded).

\section{Numerical tests}
\subsection{Drag on an ellipsoid}
In order to verify that we can accurately resolve the radius of a slender object in a fluid, we measure the drag on a slender ellipsoidal particle. The ellipsoid is positioned in the center of a periodic domain of size $[-L/2,L/2]^3$, where the limit $L \rightarrow \infty$ will be studied computationally. The half-minor axis is $a=1.33/64$ $\mu$m (chosen to be the hydrodynamic radius of the fluid grid with spacing $1/64$ $\mu$m). The half-major axis is $b=1/2$ $\mu$m. 

We treat the ellipsoid as a 1D slender body in the following way: consider the ellipsoid as a fiber with radius $R(s)$, where $0 \leq s \leq 2b$ is the arclength parameter along the major axis. Then we define the position and radius of the ellipsoid centerline as 
\begin{equation}
\label{eq:ellipseR}
\bm{X}(s) = (s-b,0,0), \qquad R(s)=\frac{\sqrt{s(2b-s)}}{2\beta}, 
\end{equation}
where $\beta = b/a$ is the aspect ratio. Substituting this form of $R(s)$ into Eq.\ \eqref{eq:Rs} and Eq.\ \eqref{eq:dragc}, 
\begin{equation}
\label{eq:RsEl}
R_c(s) = \frac{R_hR(s)}{R_h-R(s)}, \qquad \xi(s) = \frac{1}{6\pi \mu R_c(s)}. 
\end{equation}
We can therefore model the slender ellipsoid as a one-dimensional line of markers using the Stokes drag IB model. In the update equation, Eq.\ \eqref{eq:hybridupdate}, the drag coefficient $\xi$ becomes a function of position along the centerline of the ellipsoid.

So that $R(s) > 0$, we discretize the ellipsoid from $\Delta s$ to $(2b-\Delta s)$ with point spacing $\Delta s \approx R_h$. On each of the discrete points, we apply a uniform force of unit magnitude $\hat{\bm{F}}$. We then spread the force to the grid of spacing $h$. Because there cannot be any net force on a periodic domain when solving the Stokes equations, we distribute a uniform body force onto the entire fluid grid that makes the total force zero. We then solve the fluid equations, which have a nonzero velocity far away from the fiber. We denote this far field velocity by $\bm{u}^\infty= \bm{u}(0,-L/2,0)$ (i.e. we assign it to be the velocity at $x=0$, $y=-L/2$ in the plane of the fiber centerline $z=0$). Given velocity $\bm{U}(s)$ from Eq.\ \eqref{eq:hybridupdate}, the drag coefficient is computed for a given periodic domain size by
\begin{equation}
\label{eq:ibxi}
\left(\frac{F}{U}\right)_{IB} = \frac{N}{(\bar{\bm{U}}-\bm{u}^\infty) \cdot \hat{\bm{F}}}.
\end{equation}
Here $N$ is the number of points discretizing the ellipsoid, so that the total force magnitude is $N \cdot 1=N$, and $\bar{\bm{U}}$ is the average velocity computed from Eq.\ \eqref{eq:hybridupdate} for the hybrid IB method with $\xi=\xi(s)$ given in Eq.\ \eqref{eq:RsEl}.

Throughout this section, we compare our results with slender body theory, which gives
\begin{equation}
\label{eq:Lsbt}
\bm{U} = \frac{1}{8\pi \mu}\left(-\log{\left(\frac{e}{\rthr{(2\beta)^2}}\right)}\left( \bm{I}+\bm{X}_s\bm{X}_s\right) + 2\left(\bm{I}-\bm{X}_s\bm{X}_s\right)\right)\bm{f}(s), 
\end{equation}
where $\bm{X}_s=(1,0,0)$ is the filament tangent vector and $\bm{f}(s)$ is the force density on the fluid from the filament. Note that the non-local integral term associated with Eq.\ \eqref{eq:Lsbt} is zero for a straight filament with constant forcing (see \cite{ts04} for details on this). We emphasize that slender body theory is asymptotically accurate in $1/\beta$ (with error $\mathcal{O}(\beta^{-2}\log{\beta})$) for an ellipsoid with shape given by Eq.\ \eqref{eq:ellipseR} \cite{johnson80}.

\subsubsection{Drag parallel to major axis}

For motion parallel to the major axis, the analytical Oberbeck drag formula for an ellipsoid in Stokes flow is 
\begin{equation}
\label{eq:litxi}
\frac{F}{U} = 6\pi \mu a K, 
\end{equation}
where $a$ is the radius of the minor axis of the ellipsoid and $K$ is a shape factor given by
\begin{equation}
K=\frac{\frac{4}{3}\left(\beta^2-1\right)}{\frac{(2\beta^2-1)}{\sqrt{\beta^2-1}}\log{\left(\beta+\sqrt{\beta^2-1}\right)}-\beta},
\end{equation}
where again $\beta = b/a$ is the aspect ratio of the ellipsoid (ratio of major to minor axis) \cite{oberbeck1876uber, chwang1976hydromechanics}.

\begin{figure}
\centering     
\subfigure[Parallel to major axis]{\label{fig:eldrag}\includegraphics[width=70mm]{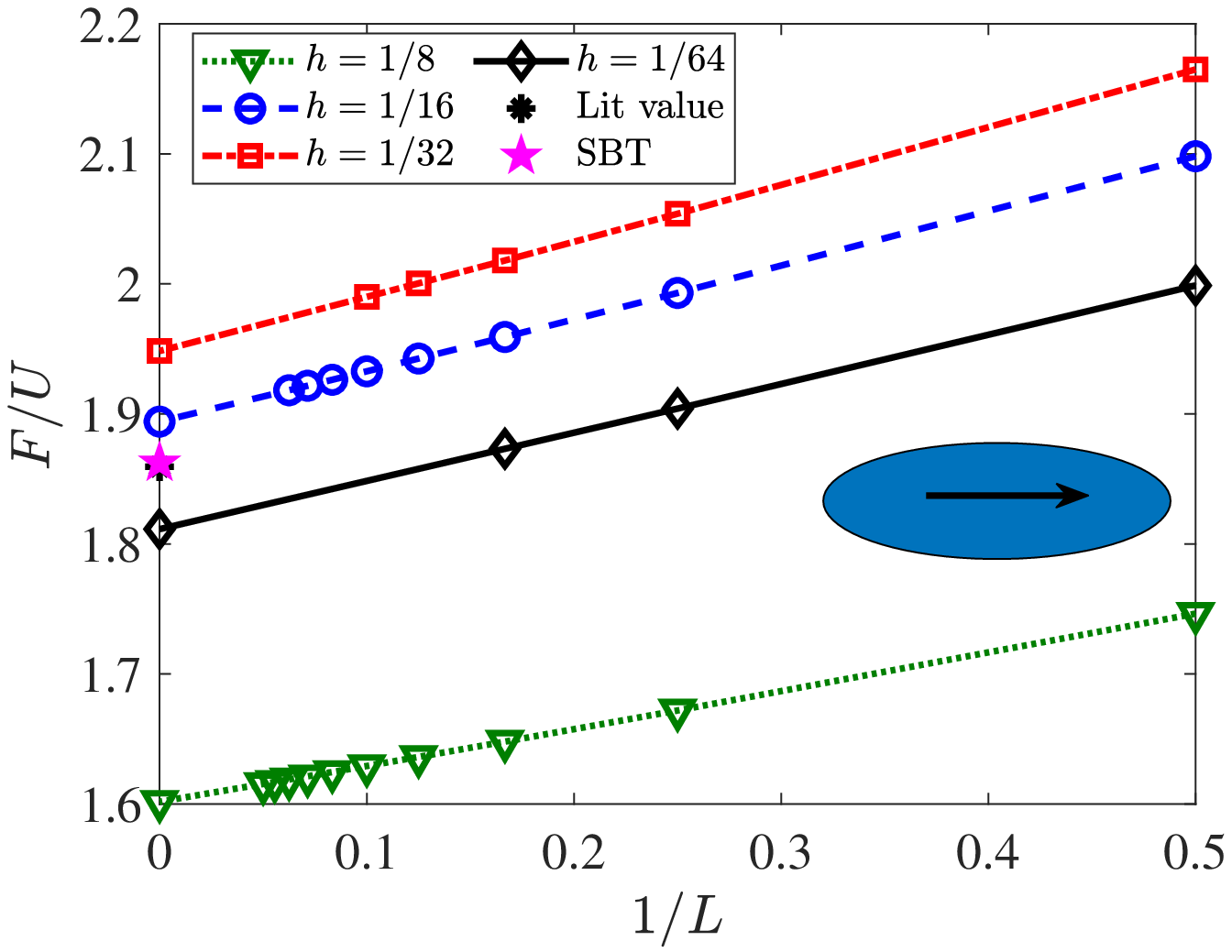}}
\subfigure[Perpendicular]{\label{fig:eldragopp}\includegraphics[width=70mm]{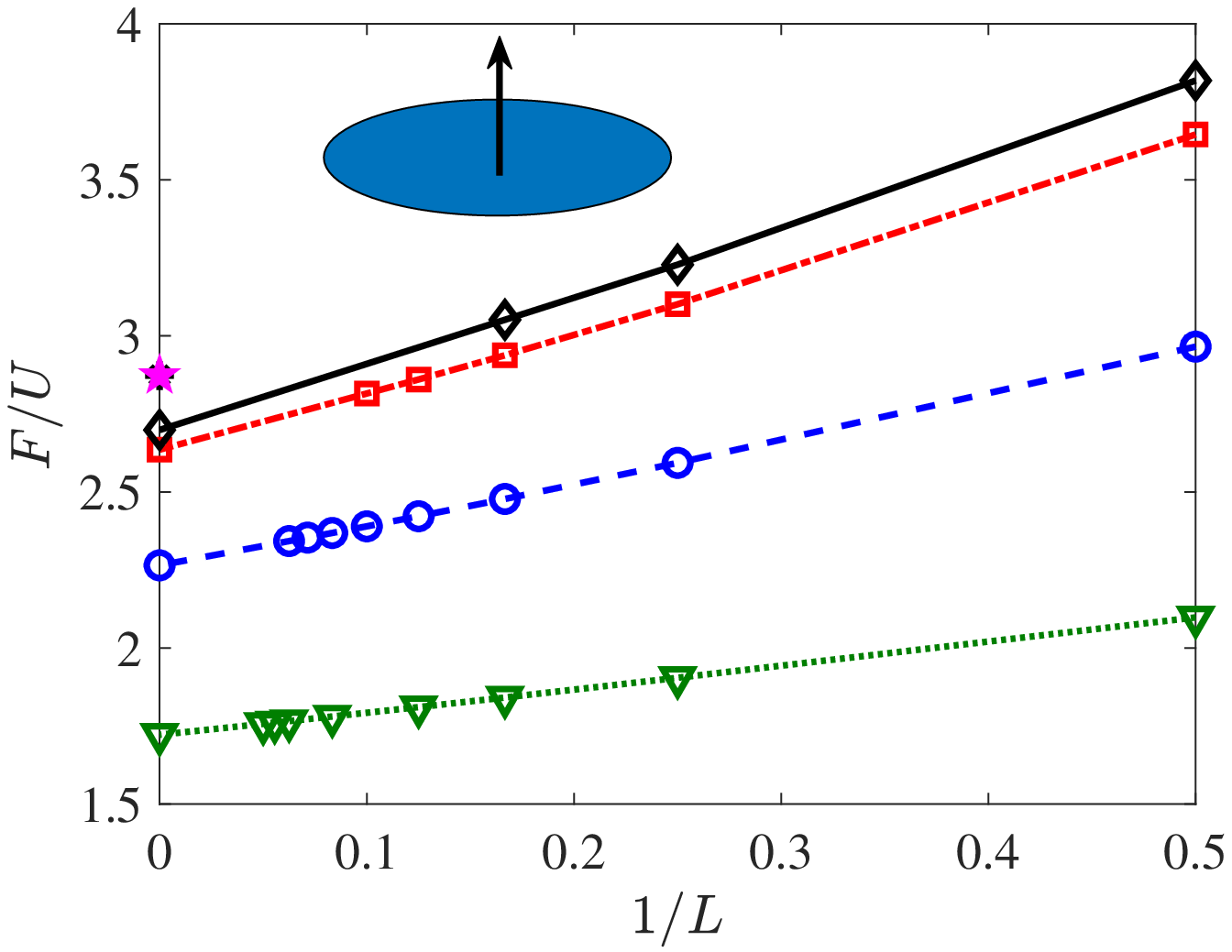}}
\caption{Computing the drag on an ellipsoid with half-minor axis $a=1.33/64$ and half-major axis $b=1/2$. $\displaystyle{\left(\frac{F}{U}\right)_{IB}}$ is shown for $h=1/8$ (green triangles), $h=1/16$ (blue circles), $h=1/32$  (red squares), and $h=1/64$ (black diamonds) in the case when the force is (a) parallel or (b) perpendicular to the major axis (see sketch on plot). We extrapolate values for finite $L$ to $L=\infty$ and compare to Eq.\ \eqref{eq:RsEl}. The values from Eq.\ \eqref{eq:RsEl} are shown as asterisks, and values from slender body theory, Eq.\ \eqref{eq:Lsbt}, are shown as a pink stars.}
\end{figure}

Fig.\ \ref{fig:eldrag} shows the computed values of $\displaystyle{\left(\frac{F}{U}\right)_{IB}}$ as a function of the inverse of the periodic domain length $L$. In order to compare our solution with the free space answer, Eq.\ \eqref{eq:litxi}, we take $L$ as large as possible until we reach a linear region in $1/L$. We then extrapolate to $L=\infty$. As shown in Fig.\ \ref{fig:eldrag}, the hybrid IB method estimates the drag coefficient to 1-2 digits of accuracy for any $h \leq 1/16$, even though the radius of the slender ellipsoid at its center is $1/64$. Note that the accuracy of the estimate \textit{does not depend strongly on the grid spacing.} This means that we have done an effective job correcting the interpolated grid velocity in Eq.\ \eqref{eq:hybridupdate}. Slender body theory, which does not involve a grid, gets three digits of accuracy on the same problem. 

\subsubsection{Drag perpendicular to major axis}

For an ellipsoid moving perpendicular to its major axis, the drag coefficient is again given by Eq.\ \eqref{eq:litxi}, but with 
\begin{equation}
K=\frac{\frac{8}{3}\left(\beta^2-1\right)}{\frac{(2\beta^2-3)}{\sqrt{\beta^2-1}}\log{\left(\beta+\sqrt{\beta^2-1}\right)}+\beta},
\end{equation}
\cite{oberbeck1876uber, chwang1976hydromechanics}.

Fig.\ \ref{fig:eldragopp} shows the computed values of $\displaystyle{\left(\frac{F}{U}\right)_{IB}}$ as a function of the inverse of the periodic domain length $L$ in this case. Although the accuracy is reduced, the hybrid IB method estimates the drag coefficient to 1 digit of accuracy for any $h \leq 1/32$. This time, the error decreases as the grid is refined, which indicates the Stokes drag terms in Eq.\ \eqref{eq:hybridupdate} at low resolutions do not accurately reflect the dynamics when the force is perpendicular to the ellipsoid center line. Equivalently, the difference in accuracy as $\bm{F}$ changes directions is indicative of the error we make in using an isotropic drag term in Eq.\ \eqref{eq:hybridupdate}. Slender body theory is again much more accurate than our IB method, this time to 4 digits. 

\subsection{Relaxing bent fibers}
\label{sec:relfibs}
\begin{figure}
\centering     
\subfigure[Set-up]{\label{fig:BendingSetUp}\includegraphics[width=70mm]{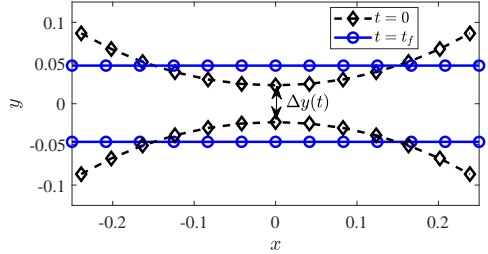}}
\subfigure[Evolution of $\Delta y(t)$]{\label{fig:bendingacc}\includegraphics[width=70mm]{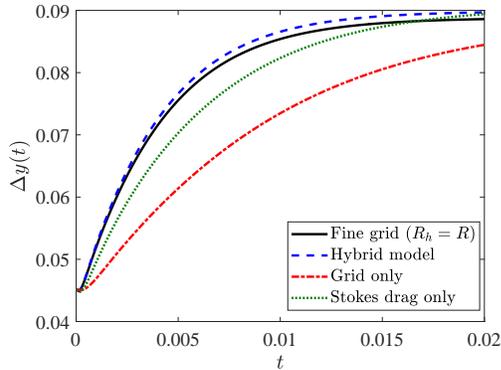}}
\caption{Numerical example demonstrating the effectiveness of the hybrid method. (a) Set up for the test. We begin with two bent fibers whose minimum separation in $y$ is $\Delta y(t)$ (black dotted lines denote initial position; diamonds are the marker point locations when $h=1/32$ $\mu$m). The initial separation $\Delta y(0)$ is $\approx 0.045$ $\mu$m so that the fibers are hydrodynamically in contact for coarser IB grids. We set the physical radius $R=0.008$ $\mu$m and relax the fibers in time; their final positions are shown as blue lines with circles. We measure the error in $\Delta y(t)$ for different grid spacings. (b) The trajectory $\Delta y(t)$ for a grid with $R_h=R$ ($h=1/168$ $\mu$m, black line) closely matches the hybrid method curve with $h=1/64$ $\mu$m (dashed blue line). The relative accuracy of 2 digits far exceeds that of the traditional IB method on a coarser grid (dashed/dotted red line) and the Stokes drag law with $R_c=R$ (dotted green line). }
\end{figure}

Having established that the hybrid method can correctly capture the radius we desire, we next seek to study its advantages in time-stepping. We thus lay out a simple numerical example to study the stability and accuracy of the hybrid method as compared to traditional IB formulations and the Stokes drag law. We consider two fibers of length $L=0.5$ $\mu$m and radius $R=0.008$ $\mu$m (aspect ratio 62.5) that are positioned in the $xy$ plane relaxing opposite each other in a domain of size $[-0.5,0.5] \times [-0.5,0.5] \times [-0.125,0.125]$ $\mu$m. The $x$ and $y$ domain sizes are chosen to avoid effects from periodicity; while the $z$ domain size is chosen to be small since we position the fibers in the $xy$ plane. Fig.\ \ref{fig:BendingSetUp} shows the initial and final position of the fibers after $t=0.02$ seconds. 

We will study the minimum separation of the fibers in time, $\Delta y(t)$, for different grid sizes and timesteps. For each meshwidth $h$, we compare the results from the coarse grid with no corrections and our method of Eq.\ \eqref{eq:hybridupdate}, using as an exact solution a grid where the fiber is totally resolved, $R_h=R$. We set the number of points on the fiber $N\approx L/R_h$, so that the points are spaced approximately one hydrodynamic radius apart \rone{(this means that the delta function support of each point modestly overlaps on the grid, which is actually recommended in the IB method)} \cite{tbring08}. For the case in which there is no grid, the points are spaced one \textit{physical} fiber diameter apart, $N \approx L/(2R)$ when $R_h=\infty$, so that the fiber is represented by a series of spheres whose edges are in contact. We study both the accuracy and stability of the method in this context.

The fibers have a discrete stretching energy
\begin{equation}
\label{eq:Es}
E_s[\bm{X}_1, \dots \bm{X}_N] = \frac{K_s}{2} \sum_{k=1}^{N-1} \left( \frac{\norm{\bm{X}_{k+1}-\bm{X}_k}}{\Delta s}-1\right)^2 \Delta s
\end{equation} 
and bending energy
\begin{equation}
\label{eq:Eb}
E_b[\bm{X}_1, \dots \bm{X}_N] = \frac{K_b}{2} \sum_{k=2}^{N-1} \frac{\norm{\bm{X}_{k+1}-2\bm{X}_k+\bm{X}_{k-1}}^2}{\Delta s^2} \, \Delta s.
\end{equation} 
In both cases, $\bm{X}_k$ denotes the $k$th point on the the fiber $\bm{X}$ and $\Delta s$ is the spacing between marker points in the reference configuration of the fibers. The force (not force density) applied by each fiber point $\bm{X}_k$ to the fluid is given by
\begin{equation}
\label{eq:force}
\bm{F}\left(\bm{X}_k\right)=\bm{F}_s\left(\bm{X}_k\right)+\bm{F}_b\left(\bm{X}_k\right)=-\frac{\partial E_s}{\partial \bm{X}_k}-\frac{\partial E_b}{\partial \bm{X}_k}
\end{equation}

We set the bending stiffness $K_b=0.25$ pN $\cdot \mu$m$^2$ (which is between the $K_b$ values of a short microtubule and actin fiber \cite{van2008microtubule, ott1993measurement}). \rone{Non-dimensionally, the timescale of fiber relaxation is $\displaystyle{t^* = \frac{\mu L_f^2 \ell^2}{K_b}}$, where $\ell$ is the approximate deviation of the fibers from their straight configuration (in this example, $\ell \approx 0.1$ and $t^* = 0.01$). To model biological filaments that are nearly inextensible, we set $K_s=100$ pN so that $\displaystyle{\frac{K_s t^*}{\mu L_f^2}\approx 10 \gg 1}$, and we find empirically that the maximum change in filament length is $\approx 1$ \%}. 

We consider three different temporal discretizations of Eq. \eqref{eq:hybridupdate}. The first is an explicit treatment; that is, 
\begin{equation}
\label{eq:expbt}
\frac{\bm{X}^{n+1}-\bm{X}^n}{\Delta t} = \bm{\mathcal{S}}^*\left(\bm{X}^n\right)\bm{u}^n+ \xi \left(\bm{F}^n_{s}+\bm{F}_b^n\right). 
\end{equation}
We next consider treating the bending force implicitly. Because bending is a linear operator on the fiber positions, we can write $\bm{F}_b = \bm{B}\bm{X}$, and therefore it is straightforward to handle the bending term implicitly. In this scenario, Eq.\ \eqref{eq:gendisc} becomes
\begin{equation}
\left(\frac{\bm{I}}{\Delta t}-\xi \bm{B}\right)\bm{X}^{n+1}=\frac{\bm{X}^n}{\Delta t}+\bm{\mathcal{S}}^*\left(\bm{X}^n\right)\bm{u}^n + \xi \bm{F}^n_{s}.
\label{eq:bendimp}
\end{equation}
We solve Eq.\ \eqref{eq:bendimp} directly (and do so for each fiber separately) because it is a sparse, banded linear system. 

The final possible temporal discretization is to treat both the tension and bending force implicitly. Because the tensile force is a nonlinear function of position, however, Newton's method is required. If we write the linearized operator
\begin{equation}
\bm{F}_{s}\left(\bm{Y} + \delta \bm{Y}\right) = \bm{F}_{s}\left(\bm{Y}\right) + \frac{\partial \bm{F}_{s}}{\partial \bm{X}}\left(\bm{Y}\right) \cdot \delta \bm{Y} + \mathcal{O}(\delta^2), 
\end{equation}
then the formula for Newton's method in this case is
\begin{align}
\label{eq:NewtonT}
& \frac{\bm{X}^{n+1, m+1} - \bm{X}^n}{\Delta t} =\bm{\mathcal{S}}^*\left(\bm{X}^n\right) \bm{u}^n\\[2 pt] \nonumber &+ \xi \left(\bm{B} \bm{X}^{n+1, m+1} + \bm{F}_{s}\left(\bm{X}^{n+1, m}\right) +\frac{\partial \bm{F}_s}{\partial \bm{X}}\left(\bm{X}^{n+1, m}\right) \cdot \left(\bm{X}^{n+1, m+1}-\bm{X}^{n+1, m}\right)\right). 
\end{align}
The calculation of the matrix $\displaystyle{ \frac{\partial \bm{F}_s}{\partial \bm{X}}}$ is given in Appendix B. We can therefore solve Eq.\ \eqref{eq:NewtonT} for $\bm{X}^{n+1, m+1}$ at each Newton iteration until $\bm{X}^{n+1}$ solves 
\begin{equation}
\label{eq:exnext}
\frac{\bm{X}^{n+1} - \bm{X}^n}{\Delta t} = \bm{S}^*(\bm{X}^n)\bm{u}^n + \xi \left(\bm{B} \bm{X}^{n+1} + \bm{F}_{s}\left(\bm{X}^{n+1}\right)\right)
\end{equation}
to some specified tolerance. We will solve Eq.\ \eqref{eq:exnext} to 6 digits of absolute accuracy.

We therefore have four possible schemes to consider: first, the explicit IB method with $\xi=0$. Next, the hybrid Stokes drag method with any of the following: explicit forcing (Eq.\ \eqref{eq:expbt}), implicit bending only (Eq.\ \eqref{eq:bendimp}), and implicit bending and stretching (Eq.\ \eqref{eq:exnext}). 

\subsubsection{Stability}
\label{sec:stability}
Table \ref{tab:twofibsstab} shows the stability of the hybrid method under the the four possible scenarios for several different grid sizes. The final grid, $h=1/168$ $\mu$m, is the grid with $R_h=R$. That is, if we wanted to use the traditional IB method to simulate the scenario in Fig.\ \ref{fig:BendingSetUp}, we would be forced to use $h=1/168$. In addition to being spatially costly, such a simulation has a severe restriction on the timestep, as the maximum stable timestep is $10^{-7}$ s for a process that takes $\mathcal{O}(10^{-2})$ seconds.

We therefore consider the hybrid method. Our first observation is that, for a given grid size, adding the Stokes drag term in an explicit fashion, Eq.\ \eqref{eq:expbt}, leads to a maximum stable timestep on the same order as the $h=1/168$ grid. This indicates that the explicit timestep is dictated by the radius of the fiber that we are attempting to resolve. Adding the Stokes drag term for a given grid size decreases the simulated fiber radius, thereby restricting the timestep even further. \rone{Equivalently, for coarser grids the Stokes drag term dominates the numerical stiffness of the problem}. 

We overcome this difficulty by treating the Stokes drag term implicitly. Indeed, for a fixed grid size, Table \ref{tab:twofibsstab} shows that we can recover the maximum stable timestep for the \textit{grid} (i.e. the same timestep as if we were to simulate a fiber with $R=R_h$ on the given grid) by treating the bending force implicitly. \rone{Sometimes, we can even improve the timestep over that needed on the grid without corrections. For example,} taking a grid size of $h=1/64$, the stable timestep with $\xi=0$ is $10^{-5}$, but adding the implicit $\xi$ term leads to a timestep a factor of 10 larger than the original grid, even though the fiber now has the same radius as on a much finer grid. \rone{This phenomenon disappears} once $R_h \rightarrow R$; as $\xi \rightarrow 0$, the contribution of the Stokes drag term in Eq.\ \eqref{eq:hybridupdate} is minimal, \rone{and the numerical stiffness of the problem is dominated by the fluid grid}. Therefore, implicit treatment of the bending force leads to little or no reduction of the timestep \rone{for finer grids} (e.g. a grid of size $h=1/128$).

\rone{With regard to implicit treatment of the stretching force, other than for the coarsest grid ($h=1/32$), we find that treating tensile forces implicitly leads to no gain in the stable timestep. This indicates that the tensile forces are \textit{not} the stiffest part of the calculation for $h < 1/32$. Indeed, when $h < 1/32$, we confirm this by decreasing $K_s$ and finding no change in the maximum stable timestep. We conclude that, for $h < 1/32$, the overall numerical stiffness is again dictated by the transmission of the bending forces through the fluid. At this point, we reach the limit of improving the maximum stable timestep via local corrections.}

\rone{In general, it seems clear from Table \ref{tab:twofibsstab} that the inclusion of the Stokes drag term has the most beneficial effects for the overall stability when $R_h$ and $R_c$ are comparable in magnitude, i.e. when the mobility is split evenly between the grid and the local corrections. }

\begin{table}
\centering
\begin{tabular}{c|c|c|c|c|c|c} 
$h$ & $R_h$ & $R_c$ & Normal IB & $\xi$, explicitly & Imp B & Imp BT \\ \hline
1/32 & 0.042 & 0.010 & $10^{-4}$ & $10^{-6}$ & $10^{-4}$ & $10^{-2}$\\
1/64 & 0.021 & 0.013 & $10^{-5}$ & $10^{-6}$ & $10^{-4}$ & $10^{-4}$\\
1/128& 0.010 & 0.035 & $10^{-6}$ & $10^{-7}$ & $10^{-6}$ & $10^{-6}$\\
1/168 & 0.008 & $\infty$ & $10^{-7}$ & - & - & - 
\end{tabular}
\caption{Stable timesteps for the hybrid method in the two fiber relaxation test, Fig.\ \ref{fig:BendingSetUp}, for varying grid sizes. ``Imp B'' denotes the treatment of bending forces implicitly via Eq.\ \eqref{eq:bendimp}, and ``Imp BT'' denotes the treatment of both bending and stretching forces implicitly via Eq.\ \eqref{eq:exnext}. \rone{Note that the relevant timescale to of fiber relaxation is $t^*\approx 0.01$ s and the physical fiber radius is $R=0.008$ $\mu$m.}}
\label{tab:twofibsstab}
\end{table}

\subsubsection{Accuracy}
So far, we have established that the stable timestep for the hybrid method is dictated by the background fluid grid size. In this section, we need to show that we can get solutions using the hybrid method over a coarser grid which are just as accurate as using the traditional IB method over a finer grid. To do this, we consider the evolution of the one-dimensional variable $\Delta y(t)$ shown in Fig.\ \ref{fig:BendingSetUp}. 

Shown in Fig.\ \ref{fig:bendingacc} is $\Delta y(t)$ over the grid with spacing $h=1/64$ using $\xi=0$ and $\xi > 0$. We see that the hybrid method solution, with $\xi > 0$, is much closer to the solution curve over the finer grid $h=1/168$ than the traditional IB method over the coarser grid ($\xi=0$). When the coarse grid is used without the hybrid method, the fibers are initially much closer together numerically since their numerical radii are larger, and so they take longer to come apart and relax to their straight configurations. But when the correction term is added, we see that we are able to almost exactly reproduce the dynamics of the traditional IB method over a fine grid. This comes at a fraction of the cost; the fine grid is over twice as fine as the $h=1/64$ grid, and the required timestep is 3 orders of magnitude smaller. Thus the ratio of costs is approximately 8000 to 1.

Finally, Fig.\ \ref{fig:bendingacc} shows that the hybrid model is also more accurate when compared to the other extreme, only using the Stokes drag model with $R_c=R$ (and markers spaced $2R$ apart; this seems to be the best choice for the marker spacing in a pure Stokes drag model). Thus by combining the IB and Stokes drag models, we are able to get more accurate results at a fraction of the cost of the traditional IB method. 

\begin{figure}
\centering     
\includegraphics[width=70mm]{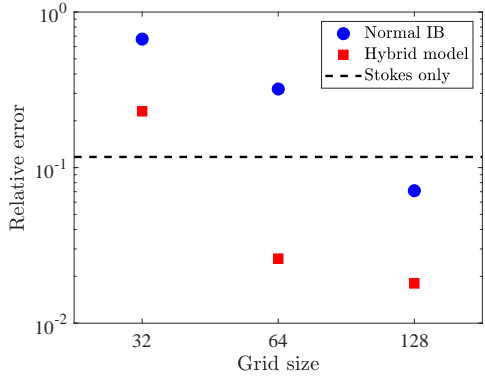}
\caption{Relative errors (maximum error in $\Delta y(t)$ divided by $\Delta y(0)$) as a function of grid size for the traditional IB method (blue circles) and hybrid IB method (red squares). Accuracy is compared to the Stokes drag model with spheres $2R$ apart (dashed black line). The hybrid model simulations are considered over the parameter set with maximum stable timestep (i.e. for $h=1/32$, implicit bending and tension, for $h=1/64$ and $h=1/128$ implicit bending only). We use as an exact solution the trajectory for a fine grid where the fiber is completely resolved ($h=1/168$).  }
\label{fig:bacc}
\end{figure}

Fig.\ \ref{fig:bacc} shows the maximum errors in $\Delta y(t)$ for several different grid spacings. The relative error is obtained by dividing the error by $\Delta y(0)$. We use as an ``exact solution'' the trajectory for a fine grid ($h=1/168$) where the fiber is completely resolved. 

We see that we can get two digits of relative accuracy on a grid ($h=1/64$) that is over 2 times as coarse as the grid required to resolve the fiber ($h=1/168$). For perspective, we notice that the accuracy of the normal IB calculation on the $h=1/64$ grid is over 10 times worse than the hybrid method. This trend is the same for grids of size $h=1/32$ and $h=1/128$, as the hybrid method is again more accurate than the traditional IB method in all cases. Finally, we observe that the hybrid method is much more accurate than the Stokes drag law alone (dotted black line), which has been integrated explicitly with $\Delta t =10^{-8}$ for maximum accuracy. 

In sum, this section has shown that the hybrid method can approximate a solution on a finer grid to within 1-2 digits of relative accuracy. This comes at a fraction of the cost of actually doing the finer-grid computation, as the timestep is dictated by the \textit{coarser} grid. For biological and engineering applications, 1-2 digits of relative accuracy can be enough to identify trends within experimental noise thresholds. We also established that the hybrid method is most effective for stability and accuracy when the mobility is approximately split evenly between the IB calculation and Stokes drag law. Although the method is not restricted to this case, it is most effective therein. 

\section{Applications}
\subsection{Fiber in shear flow \label{sec:fibshear}}
Having \rone{carried out} a detailed numerical example to show the advantages of the method, we next apply it to a \rone{ benchmark test case} to show we can reproduce existing results. We consider a single fiber in a shear flow. Previous studies, reviewed in \cite{du2019dynamics}, have shown that there are three characteristic modes by which fibers deform in a shear flow: tumbling, in which the fiber remains straight and tumbles end-over-end, buckling, in which the fiber initially appears to be tumbling before buckling in the middle, and U-snaking, in which one edge of the fiber leads the motion by folding over the other edge. The transition between modes is governed by a single dimensionless number, defined in \cite{quennouz2015transport} to be the elasto-viscous number, 
\begin{equation}
\tilde{\eta} = \frac{\mu L_f^4 \dot{\gamma}}{K_b},
\end{equation}
where $\mu$ is the fluid viscosity, $L_f$ is the fiber length, $\dot{\gamma}$ is the shear rate of the fluid, and $K_b$ is the bending stiffness of the fiber. In \cite{liu2018morphological}, the range of $\tilde{\eta}$ for each mode is determined from experiments and simulations to be approximately $\tilde{\eta} \in [0, 300]$ for tumbling, $\tilde{\eta} \in [300, 2000]$ for buckling, and $\tilde{\eta} \in [2000, 20000]$ for U-snaking. 

In this study, we choose three representative elasto-viscous numbers, $\tilde{\eta}=150, 450$, and 7500. Our goal is to show that we can reproduce tumbling, buckling, and U-snaking, respectively. 

We initialize a fiber of length $L_f=1$ $\mu$m inside of an Eulerian domain of size $[-1,1] \times [-1,1] \times [-1/8,1/8]$ $\mu$m. The Eulerian domain has $h=1/32$ $\mu$m, so that it is $64 \times 64 \times 8$ grid points. We set the radius of the fiber to be $R=R_h/2$, where $R_h=1.33/32$ $\mu$m is the approximate hydrodynamic radius of the grid. We have done this so the contribution from the traditional IB terms and Stokes drag terms in Eq.\ \eqref{eq:hybridupdate} are equal, i.e. $R_c=R_h=2R$. 

We fix $\mu=1$ Pa$\cdot$s and $\dot{\gamma}=3$ s$^{-1}$, so that $\tilde{\eta}$ is completely determined from $K_b$. In particular, we have $K_b=0.02, 0.0067, 4 \times 10^{-4}$ pN $\cdot$ $\mu$m$^2$ for the values $\tilde{\eta}=150, 450, 7500$ that we have chosen. We fix $K_s=100$ pN to preserve fiber inextensibility (throughout a 10 second simulation in all cases, the change in fiber length was $\leq 0.5$\%).

We initialize the fiber to be straight with a small perturbation to break the symmetry,
\begin{equation}
\bm{X}(s) = \begin{pmatrix} s-\frac{L}{2}\\[2 pt] 0\\[2 pt] 0 \end{pmatrix}+0.1 \begin{pmatrix}0\\[2 pt] e^{-s}-\langle e^{-s} \rangle \\[2 pt] 0\end{pmatrix}. 
\end{equation} 
Here $s \in [0,L]$ is the arclength parameter, and the second term is such that the fiber is centered on 0 in the $y$ direction (we subtract the mean of $e^{-s}$). We implement a shear flow in an immersed boundary framework by splitting the Stokes equations into two parts. Specifically, we define
\begin{equation}
\label{eq:ush}
\bm{u}_{sh} = \dot{\gamma} \begin{pmatrix} y \\ 0 \\ 0 \end{pmatrix}. 
\end{equation}
and write 
\begin{equation}
\bm{u} = \bm{u}_{sh} + \bm{u}_p, 
\end{equation}
where $\bm{u}_p$ \rthr{is the velocity field obtained by solving the Stokes equations as driven by the fiber force density on a periodic domain. In this formulation, $\bm{u}_{sh}$ satisfies the free space Stokes equations by construction, and the periodic domain is assumed to be large enough relative to the fiber so that periodic effects are negligible and we are approximately solving the free space problem of a single fiber in a shear flow}.
We therefore add the linear shear $\bm{u}_{sh}$ to the formulation of Eq.\ \eqref{eq:gendisc} by writing
\begin{equation}
\label{eq:saoseqn}
\frac{\bm{X}^{n+1}-\bm{X}^n}{\Delta t} = \bm{\mathcal{S}}^*(\bm{X}^n)\bm{u}^n + \bm{u}_{sh}^n+\xi \bm{F}(\bm{X}^{n+1}),
\end{equation}
and solving Eq.\ \eqref{eq:saoseqn} on each fiber. Since $\bm{u}_{sh}$ is a zeroth order term (i.e. it does not involve derivatives of $\bm{X}$), it adds virtually no additional stiffness to the problem and can be treated explicitly. Similar to \cite{ts04}, we simulate only a half cycle of fiber motion, as in the absence of thermal fluctuations it is impossible for the fiber to leave its flat steady state after the half cycle has been completed. 

We note the convenience of the hybrid method in this situation. We simulate 10 seconds of fiber deformation. Using the \textit{hybrid method} with implicit tension and bending, Eq.\ \eqref{eq:exnext}, the maximum stable timestep is $\Delta t =10^{-3}$ s, so that each simulation takes only seconds to run on a laptop. Meanwhile, if we performed these simulations by the traditional IB method, the grid would be $128 \times 128 \times 16$ to capture the correct fiber radius, and the maximum stable timestep would drop to $10^{-5}$ s. Each simulation is thus $800$ times less expensive than with the traditional IB method. 

\begin{figure}
\centering     
\subfigure[$\tilde{\eta}=150$, Tumbling]{\label{fig:Tumbling}\includegraphics[width=70mm]{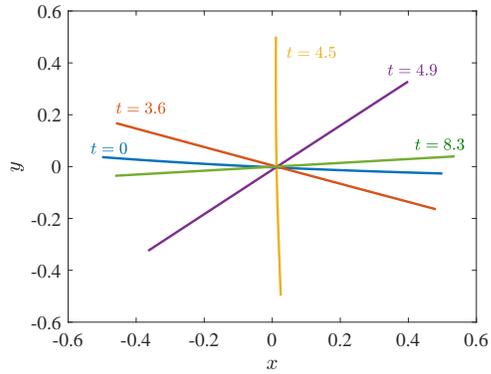}}
\subfigure[$\tilde{\eta}=450$, Buckling]{\label{fig:Buckling}\includegraphics[width=70mm]{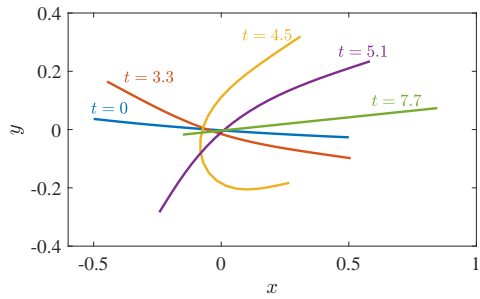}}
\subfigure[$\tilde{\eta}=7500$, U-snaking]{\label{fig:USnake}\includegraphics[width=70mm]{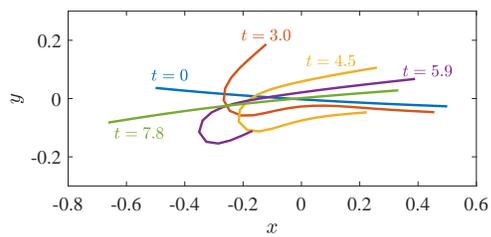}}
\caption{Fibers in shear flow with varying elasto-viscous number $\tilde{\eta}$. (a) $\tilde{\eta}=150$ shows end-over-end tumbling behavior. (b) $\tilde{\eta}=450$ shows buckling in a C shape, (c) $\tilde{\eta}=7500$ demonstrates U-snaking. }
\label{fig:fibshear}
\end{figure}

We see in Fig.\ \ref{fig:fibshear} that the hybrid IB method is able to reproduce behavior from the phase diagram for varying elasto-viscous number. Fig.\ \ref{fig:Tumbling} shows that when the bending stiffness is large enough ($\tilde{\eta}=150$), the fiber tumbles end-over-end without buckling. When the bending stiffness is reduced  ($\tilde{\eta}=450$), Fig.\ \ref{fig:Buckling} shows that the fiber is deformed into a C-shape as it buckles under the shear flow. Finally, for very low bending stiffness  ($\tilde{\eta}=7500$), the initially raised edge of the fiber at $x=-0.5$ is turned almost independently from the rest of the fiber. After the leading edge is sufficiently bent, the lagging edge follows, and the fiber follows a U-snaking pattern. Thus the hybrid method is able to reproduce the behavior of a single fiber in a shear flow for a range of elasto-viscous numbers. 

\subsection{Viscosity of an aligned fiber suspension}
One of the advantages of the hybrid method is that it scales linearly with the number of fibers. In particular, because the linear solve in Eq.\ \eqref{eq:gendisc} is done on a per fiber basis, it is computationally much easier to simulate a suspension of a large number of fibers than in the case in which the non-local hydrodynamic interactions are treated implicitly. With this in mind, we focus in this section on the simulation of a suspension of high-aspect-ratio fibers. In particular, our goal is to measure the viscosity of a fiber suspension and compare the results to those of \cite{kochsumfig, fibexps, mackfibs}, which give a summary of results for fibers aligned by a linear shear flow for a variety of volume fractions. 

\begin{figure}
\centering     
\includegraphics[width=70mm]{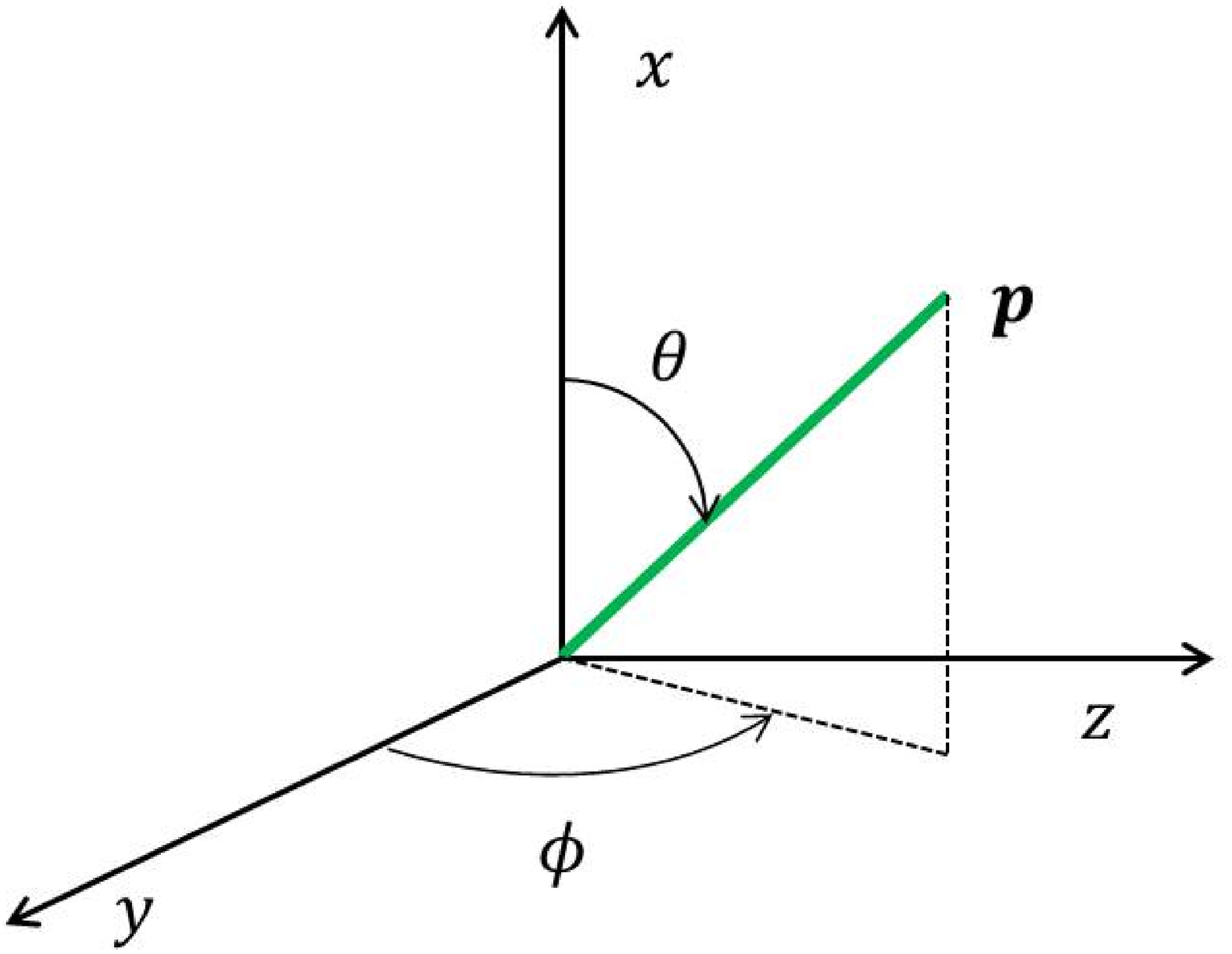}
\caption{Spherical coordinate system for the fiber orientation distribution, $\Omega$. $x$ is the flow direction, $y$ is the gradient direction, and $z$ is the vorticity direction.}
\label{fig:sphcoords}
\end{figure}

Let $n$ be the number of fibers per unit volume and $L_f$ and $d$ be the length and diameter of a fiber, respectively. Then, in the dilute ($nL_f^3 \ll 1$) and semi-dilute ($nL_f^2 d \ll 1$) regimes, the distribution function for the orientation vector (on the sphere) of fibers subject to a simple shear flow has been derived analytically using slender body theory and taking into account one- and two-fiber interactions \cite{kochsumfig, kochdist1, kochdist2}. Consider the coordinate system in Fig.\ \ref{fig:sphcoords}, where $x$ is the flow direction, $y$ is the gradient direction, and $z$ is the vorticity direction. Then the distribution of orientations $\bm{p}$ on the unit sphere is given by
\begin{equation}
\label{eq:omega}
\Omega(\bm{p}) = \frac{R}{\pi r_e \left(4R\left(r_e^{-2}\cos^2{\theta} + \sin^2{\theta}\cos^2{\phi}\right)+\sin^2{\theta} \sin^2{\phi}\right)^{3/2}}. 
\end{equation}
Here $R$ is a constant (generally determined from a best fit to a combination of experimental and theoretical data) that has to do with the anisoropic diffusivity of the orientation vector on the unit sphere. We set $R=3$, which is the experimentally measured value for a semi-dilute suspension with an aspect ratio $r_p=L_f/d \approx 33$ \cite[Table~1]{kochdist2}. Derivations of the orientation functions are generally done using Jeffrey orbits, which are based on ellipsoidal fibers. Therefore, in Eq.\ \eqref{eq:omega}, the parameter $r_e \approx 0.7r_p$ is an adjustment to the aspect ratio to approximate the orientation distribution for \textit{cylindrical} fibers \cite{kochsumfig, kochdist1, kochdist2}.

The goal here is to show that our numerical method can reproduce existing measurements on the viscosity of a suspension with orientation distribution given in Eq.\ \eqref{eq:omega}. These measurements come from both experiments \cite{fibexps} and direct numerical simulations \cite{kochsumfig}. Since many of the available measurements are for $r_p \approx 33$, we also take $r_p=33$ with $L_f=0.5$ $\mu$m. The resulting fiber radius is $7.5 \times 10^{-3}$ $\mu$m. We consider the fibers on an Eulerian grid with spacing $h=1/64$ $\mu$m, and so $R_{h}=0.021$ $\mu$m and $R_c=0.012$ $\mu$m in Eq.\ \eqref{eq:Rs}. The fibers are initialized as straight fibers within the periodic domain, and are assumed to be nearly inextensible, so that $K_s = 100$ pN with $K_b = 0.25$ pN $\cdot \mu$m$^2$. In Eq.\ \eqref{eq:gendisc}, we treat bending implicitly via Eq.\ \eqref{eq:bendimp} (the timestep is $\Delta t =10^{-6}$). We treat stretching explicitly, however, since we find that there is not much benefit in treating it implicitly. 

We initialize anywhere from 40 ($nL_f^3=5$) to 640 ($nL_f^3=80$) fibers inside of the periodic Eulerian domain $[-0.5,0.5]^3$ in the following way. For fiber $i$, we begin by randomly choosing a start location $\bm{s}_i$ from a uniform distribution over the periodic Eulerian domain $[-0.5,0.5]^3$. We next choose the fiber orientation vector $\bm{p}$ by sampling from $\Omega(\bm{p})$ in Eq.\ \eqref{eq:omega}. We do this by rejection. That is, we compute $\theta = \pi \times \text{rand}$, $\phi = 2\pi \times \text{rand}$, $q = \Omega_{max} \times \text{rand}$ and accept a possible fiber orientation vector $\bm{p}_i$ if the randomly chosen $q < \Omega(\theta,\phi) \sin{\theta}$. Note the Jacobian factor of $\sin{\theta}$ which comes from the fact that $\Omega$ is a distribution over the unit sphere. Once a fiber orientation vector $\bm{p}_i$ and start location $\bm{s}_i$ have been chosen independently, we set the endpoint of the fiber to $\bm{e}_i=\bm{s}_i+L_f\bm{p}_i$.

Once the fibers are initialized, our goal is to measure the viscosity of the suspension. To do this, we apply a force density on the fluid of the form
\begin{equation}
\bm{f}(y) = f_0 \begin{pmatrix} \sin{\left(\frac{2\pi y}{L}\right)} \\ 0\\ 0\end{pmatrix}:=\begin{pmatrix}f(y) \\ 0\\ 0\end{pmatrix}.
\end{equation}
Under this forcing, Eq.\ \eqref{eq:Stokeseqnsone} reduces to 
\begin{equation}
\label{eq:sred}
0 = \mu \frac{\partial^2 u}{\partial y^2} + f(y).
\end{equation}
where $u$ is the velocity in the $x$ direction. Multiplying Eq.\ \eqref{eq:sred} by $f(y)$ and integrating over $y$ gives 
\begin{equation}
\mu(x,z) = \frac{f_0 L^3}{8\pi^2 \int_{-L/2}^{L/2} u(x,y,z) \sin{\left(\frac{2\pi y}{L}\right)}} \, dy
\end{equation}
as the viscosity for a fixed $(x,z)$. We set $f_0=0.1$ pN/$\mu$m$^3$ and average $\mu(x,z)$ over all $(x,z)$ to obtain a total system viscosity $\mu$. 

We set the viscosity of the underlying fluid to be $\mu_f=1 $ Pa $\cdot$ s and compute the additional viscosity provided by the fibers from
\begin{equation}
\label{eq:mueff}
\mu_{\text{eff}} =\frac{\mu}{\mu_f}-1. 
\end{equation}
At $t=0$, the fibers are not stressed, and so $\mu_{\text{eff}}=0$. Our tests showed $c(1-e^{-kt})$ type growth to the steady state value of $\mu_{\text{eff}}$. We run the simulations until our measurements for $\mu_{\text{eff}}$ (taken every 10 timesteps) are constant to 5 digits; we then assign the steady state $\mu_{\text{eff}}$ to be the final value of $\mu_{\text{eff}}$. Since we use $f_0=0.1$ pN/$\mu$m$^3$, displacements of the fibers are small (at most $\mathcal{O}(10^{-4})$ $\mu$m over the $\mathcal{O}(10^{-2})$ s timescales that we consider), and the fibers can be assumed rigid, which is the same assumption as that made when deriving the distribution $\Omega$ in Eq.\ \eqref{eq:omega} \cite{kochdist1, kochdist2}. Another way to state this assumption is that the strain in the fibers is negligible, even though they are developing stress which leads to a steady-state viscosity. 

\begin{figure}
\centering     
\includegraphics[width=70mm]{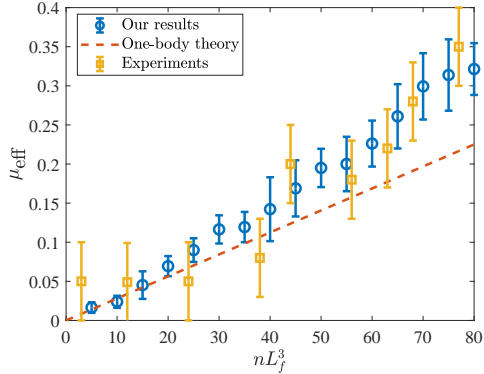}
\caption{Additional viscosity of the fiber suspension, measured using Eq.\ \eqref{eq:mueff}, as a function of $nL_f^3$. We show the average over 10 trials as blue circles, with error bars indicating the standard deviation. Also shown is a one body hydrodynamic theory \cite{kochsumfig} (dashed red line) and a set of experimental results \cite{fibexps} (yellow squares).}
\label{fig:muext}
\end{figure}

Fig.\ \ref{fig:muext} shows the measurements for effective viscosity as a function of $nL_f^3$, averaged over 10 trials (error bars show the standard deviation in the effective viscosity). We note the outstanding agreement with existing theoretical and experimental results. Fig.\ \ref{fig:muext} shows the values obtained with a one-body slender body theory \cite{kochsumfig} and experimentally \cite{fibexps}, both of which agree with our measured values. In addition, both our data and the experimental data deviate from the one body theory for large $nL_f^3$, as expected intuitively. Note also how we are able to match (within error) the large experimental viscosities for $nL_f^3 \geq 60$. For more data, see \cite[Figure~2]{kochsumfig} and \cite[Figure~10]{mackfibs}. 

\section{Conclusions}
Classical immersed boundary simulations with high aspect ratio objects can be costly because of the need to resolve competing length scales. Since the effective radius of a fiber is related to the width of the underlying fluid grid, a slender fiber will require many grid points along its length, and therefore many markers in its Lagrangian discretization. 

Here we take a step in remedying this problem by partially decoupling the movement of the fiber from the velocity on the grid. We split the velocity into two components: a drag component proportional to the force, and a classical IB component coming from a coarser grid. The drag coefficient is determined uniquely as a function of the grid spacing and physical fiber radius (i.e. the real fiber radius). Adding a drag term also allows for an explicit-implicit temporal integration scheme, in which the velocity from the grid is treated explicitly and (part of) the drag velocity is treated implicitly. This reduces the overall stiffness of IB calculations because (a) simulations can be run on coarser grids with fewer markers for the body, and (b) high order bending forces can be treated implicitly at little additional computational cost. \rtwo{The method presented here is linear in the number of particles and is effectively an inexact Ewald splitting procedure for the IB method, where the far field hydrodynamics are done on a coarse grid with a ``near field'' correction that is entirely local. While the splitting is inexact, the method presented here is an inexpensive way to obtain 1-2 digits of accuracy and can be easily implemented in existing IB solvers. }

We give several examples to show how our method improves the stability of the IB method while giving similar numerical results. We show that we can compute the drag on an \textit{ellipsoid} to 1-2 digits of accuracy simply by using a 1D line of IB markers and changing the Stokes drag coefficient in a position-dependent manner to capture the local radius of an ellipsoid. This is in contrast to a traditional IB method, where an ellipsoid must be simulated by positioning markers around its perimeter (otherwise, a line of markers will behave like a cylinder). We compare our results to slender body theory, finding the latter to be much more accurate. Simulations that only involve a small number of flexible fibers would therefore be more accurately done with slender body theory (assuming the singular integrals therein can be computed accurately), but our goal is to reach large numbers of fibers, where the quadratic-complexity interaction terms in slender body theory become significant. Furthermore, an immersed boundary framework such as this one can be used to couple slender bodies to non-slender ones on a coarse grid. 

We next use a simple example of initially bent fibers to show that the hybrid method can give improved stability over the classical IB method. In particular, on a grid twice as coarse as necessary, we obtain 2 digits of accuracy with a timestep 100 times larger than that required over the finer grid. Finally, we verify our method by applying it to two test cases: a single fiber in shear flow and a suspension of hundreds of fibers. In the latter case, we measure the viscosity of a suspension of high-aspect-ratio fibers and obtain results that agree with existing theory and experiments \cite{kochdist1, kochdist2, kochsumfig, fibexps, mackfibs}. Our simulations involve as many as 640 fibers and have a computational cost that is linear in the number of fibers.  

Note however that the key component of achieving subgrid resolution in the present method is the partial decoupling of the structure velocity from the velocity field derived from the background fluid grid.  This unfortunately erases one of the key advantages of the IB method. In a standard IB formulation, fibers that come into close contact cannot pass through each other as they are advected by a continuous velocity field obtained by interpolation from the background fluid grid. Removing this property can therefore allow fibers to run into and pass through each other. This could be remedied, if necessary, by using a repulsive force or potential when the fibers come close to each other.

\appendix
\section{Measuring the hydrodynamic radius}
In this section, we repeat the experiment of \cite{tbring08} to measure the hydrodynamic radius of our IB solver. We take an Eulerian domain of size $[-0.5,0.5]^3$ and grid size $h=1/128$ with viscosity $\mu=1$ and initialize two points inside of the domain at random locations (in the test here these two points are $(-0.31, 0.14, 0.04)$ and $(-0.46, -0.22, 0.20)$). At each timestep, we apply a random force of magnitude from 0 to 1 in each direction on the first particle, The force on the second particle is automatically the negative of the first to satisfy the constraint of no net force in periodic Stokes flow. We then use the IB method to determine the velocity of each point. Using Eq.\ \eqref{eq:sdrag}, we determine the hydrodynamic radius in each direction from 
\begin{equation}
\label{eq:hydroR}
R_h = \frac{F_{ij}}{6 \pi \mu U_{ij}}.
\end{equation}
Here the index $i$ indicates the direction $i=1,2, 3$ and $j=1,2$ indicates the particle number. Because there are 6 index combinations in Eq.\ \eqref{eq:hydroR}, we can get 6 measurements of the hydrodynamic radius at each timestep. We run this for 500 timesteps with $\Delta t = 10^{-3}$, so that by the end of 500 timesteps the points have moved approximately 2 domain sizes. We have 3000 values for $R_h$ and take the mean, finding $R_h\approx1.34h$. 
Other grid sizes give similar but distinct values; for this paper we always set $R_h=1.33h$, but we encourage this simple measurement to be done prior to implementing the hybrid method, as it can change based on the fluid solver and delta function kernel. Our value is in agreement with previous results for the 4 point delta function and a spectral fluid solver \cite{tbring08}.  

\section{Differentiating the elastic force}
Differentiating the elastic force is a non-trivial calculation, and so we present it here. Let $\bm{X}_{k-1}$, $\bm{X}_ k$, and $\bm{X}_{k+1}$ be the coordinates of a set of three points along the fiber. Let $R_{k+1} = \norm{\bm{X}_{k+1}-\bm{X}_k}$. Then the elastic force at node $k$ is given by
\begin{equation}
\label{eq:elforce}
\bm{F}_k = T\left(R_{k+1}\right) \left(\frac{\bm{X}_{k+1}-\bm{X}_k}{R_{k+1}}\right) - T\left(R_k\right)  \left(\frac{\bm{X}_{k}-\bm{X}_{k-1}}{R_{k}}\right),
\end{equation}
where as usual $T(R)=K_s(R-1)$ is the fiber tension. From this it is easy to obtain the derivative of the force with respect to $\bm{X}_{k+1}$. This $3 \times 3$ matrix is
\begin{equation}
\label{eq:dfdxlong}
\frac{\partial \bm{F}_{k, \alpha}}{\partial \bm{X}_{k+1, \beta}} = \left(T'(R_{k+1})-\frac{T(R_{k+1})}{R_{k+1}}\right)\frac{\bm{X}_{k+1,\alpha}-\bm{X}_{k, \alpha}}{R_{k+1}} \frac{\partial R_{k+1}}{\partial \bm{X}_{k+1, \beta}} + \delta_{\alpha \beta} T(R_{k+1}) \frac{1}{R_{k+1}}.
\end{equation}
Let $\displaystyle \bm{\tau}_{k+1}=\frac{\bm{X}_{k+1}-\bm{X}_k}{R_{k+1}}$. Then Eq.\ \eqref{eq:dfdxlong} can be rewritten as
\begin{align}
\frac{\partial \bm{F}_k }{\partial \bm{X}_{k+1}} & = \left(T'(R_{k+1}) - \frac{T(R_{k+1})}{R_{k+1}}\right) \bm{\tau}_{k+1}\bm{\tau}_{k+1}^* + \frac{T(R_{k+1})}{R_{k+1}}\bm{I}\\[2 pt] & = T'(R_{k+1})\bm{\tau}_{k+1}\bm{\tau}_{k+1}^* + \frac{T(R_{k+1})}{R_{k+1}}\left(\bm{I}-\bm{\tau}_{k+1}\bm{\tau}_{k+1}^*\right), 
\end{align}
where $\bm{\tau}\bm{\tau}^*$ denotes an outer product. The computation of $\displaystyle \frac{\partial \bm{F}_k}{\partial \bm{X}_{k-1}}$ is similar. $\displaystyle \frac{\partial \bm{F}_k}{\partial \bm{X}_{k}}$ can then be computed indirectly by observing that $\displaystyle \sum_j \bm{F}_j = \bm{0}$. Therefore, 
\begin{equation}
\frac{\partial \bm{F}_k}{\partial \bm{X}_k} = -\sum_{j \neq k} \frac{\partial \bm{F}_j}{\partial \bm{X}_k}. 
\end{equation}
The latter sum has at most two terms. 

\section*{Acknowledgments}
Ondrej Maxian is supported by the NSF Graduate Research Fellowship \#DGE1342536 and the Henry MacCracken fellowship. 

\bibliographystyle{siamplain}
\bibliography{SIAM_Maxian_Master}
\end{document}